\documentclass[12pt]{amsart}
\usepackage{amsmath, geometry, amssymb} 
\numberwithin{equation}{section}
\newtheorem{theorem}{Theorem}[section]

\newtheorem{corollary}{Corollary}[section]

\newtheorem{definition}{Definition}[section]

\newtheorem{example}{Example}[section]

\newtheorem{remark}{Remark}[section]

\title{Jacobi sums and correlations of Sidelnikov sequences}
\author{Ay\c{s}e Alaca and Goldwyn Millar}
\date{} 

\begin{document}

\maketitle

\begin{abstract} We consider the problem of determining the cross-correlation values of the sequences in the families comprised of constant multiples of $M$-ary Sidelnikov sequences over $\mathbb{F}_q$, where $q$ is a power of an odd prime $p$. We show that the cross-correlation values of pairs of sequences from such a family can be expressed in terms of certain Jacobi sums. This insight facilitates the computation of the cross-correlation values of these sequence pairs so long as $\phi(M)^{\phi(M)} \leq q.$ We are also able to use our Jacobi sum expression to deduce explicit formulae for the cross-correlation distribution of a family of this type in the special case that there exists an integer $x$ such that $p^x \equiv -1 \pmod{M}.$
\end{abstract}

\section{Introduction}

A sequence $\mathbf{a} = a_0 a_1 a_2 \ldots$  of elements from the ring $\mathbb{Z}/M\mathbb{Z},$ where $M$ is a positive integer, is called an $M$-ary sequence.
A sequence $\mathbf{a}$ is said to be \emph{periodic}
if there is an integer $v>0$ such that  $a_i = a_{v+i}$ for all integers $i\geq 0$.
If $v$ is the smallest such integer, then we say that $\mathbf{a}$ has period $v$. For the rest of this section, we assume that $\mathbf{a}$ is a periodic $M$-ary sequence of period $v$. 

Let $\mathbf{b} = b_0 b_1 b_2 \ldots$ be another periodic $M$-ary sequence of period $v.$ The (periodic) correlation $\mathcal{C}_{\mathbf{a},\mathbf{b}}$ of $\mathbf{a}$ and $\mathbf{b}$ is defined as follows: for each nonnegative integer $\tau$,
\begin{eqnarray}
\label{Aeqn1}
\mathcal{C}_{\mathbf{a},\mathbf{b}}(\tau) := \sum_{t = 0}^{v-1} \text{exp}\Big({\frac{2\pi i(a_t - b_{t + \tau})}{M}}\Big),
\end{eqnarray}
where the terms of $\mathbf{a}$ and $\mathbf{b}$ appearing in the exponents are interpreted as integers. 

The function $\mathcal{C}_{\mathbf{a},\mathbf{a}}$ is called the \emph{autocorrelation} of $\mathbf{a},$ and the values $\mathcal{C}_{\mathbf{a},\mathbf{a}}(\tau)$ for $1\leq \tau \leq v-1$ are called the \emph{out-of-phase} autocorrelation values of $\mathbf{a}.$ We say that $\mathbf{a}$ has \emph{low out-of-phase autocorrelation} if its out-of-phase autocorrelation values are small compared to $v.$ Low out-of-phase autocorrelation is one of the criteria for a periodic sequence to be suitable for use as a key sequence in a stream-cipher cryptosystem (see the discussion of Golomb's postulate R3 in \cite[Section 5.1]{G1} and \cite{W1}). Indeed, for a sequence to be useful for this purpose, its out-of-phase autocorrelation values should be close to zero.

If there exists an integer $\ell$ such that for each positive integer $i,$ $a_i = b_{i+\ell},$ then we say that $\mathbf{a}$ and $\mathbf{b}$ are \emph{shift-equivalent} and that $\mathbf{a}$ and $\mathbf{b}$ are \emph{shifts} of one another; in this case, $\mathcal{C}_{\mathbf{a},\mathbf{b}}$ can be obtained from 
$\mathcal{C}_{\mathbf{a},\mathbf{a}}$ by
\[\mathcal{C}_{\mathbf{a},\mathbf{b}}(\tau) = \begin{cases} \mathcal{C}_{\mathbf{a},\mathbf{a}}(\tau-\ell) & \text{ if } \tau\geq \ell, \\ \mathcal{C}_{\mathbf{a},\mathbf{a}}(\tau+v-\ell)   & \text{ if } \tau<\ell . \end{cases}\] 

If $\mathbf{a}$ and $\mathbf{b}$ are shift-inequivalent, then we say that $\mathcal{C}_{\mathbf{a},\mathbf{b}}$ is the \emph{cross-correlation} of $\mathbf{a}$ and $\mathbf{b}.$ Furthermore, we consider a family $\mathcal{F}$ of shift-inequivalent sequences to have \emph{low cross-correlation} if for any pair of sequences $\mathbf{c},\mathbf{d}$ in $\mathcal{F}$ and for any $\tau,$ $\mathcal{C}_{\mathbf{c},\mathbf{d}}(\tau)$ is small compared to $v.$ The following result, which is due to Welch \cite{W2}, gives a lower bound on the cross-correlation and out-of-phase autocorrelation values of the sequences in a family of a given size.
\begin{theorem} \cite[Theorem (periodic correlation)]{W2} \label{welch} Let $S$ be a family of $s$ shift-inequivalent $M$-ary sequences of period $\ell.$ Let $\mathcal{C}_{\text{max}}$ denote the maximum of the cross-correlation and out-of-phase autocorrelation values of the sequences in $S.$ Then 
\[\mathcal{C}_{\text{max}} \geq\sqrt{ \frac{\ell^2(s-1)}{s\ell-1}}.\]
\end{theorem}

\begin{remark}
Following \cite{G1}, we say that $\mathcal{F}$ has low cross-correlation if for any pair of sequences $\mathbf{c},\mathbf{d}$ from $\mathcal{F},$ $\mathcal{C}_{\mathbf{c},\mathbf{d}}$ outputs only values less than or equal to $\delta\sqrt{v} + \epsilon,$ for some small integers $\delta$ and $\epsilon$.
\end{remark}
 
Families of shift-inequivalent sequences with low cross-correlation (and such that each member has low out-of-phase autocorrelation) have found use in the design of code division multiple access (CDMA) radio communications systems (see, for instance, \cite{G1} or \cite{G3}). For CDMA applications, one would like such families to be as large as possible (i.e. to include as many sequences as possible). One may also desire that sequences in these families have certain additional properties (such as cryptographic strength). 

One indicator of cryptographic strength (besides low out-of-phase autocorrelation) is the \emph{balance property}. We say that the sequence $\mathbf{a}$ is \emph{balanced} if in a given period of $\mathbf{a}$ (i.e. in a given list of $v$ consecutive elements of $\mathbf{a}$) each element of $\mathbb{Z}/M\mathbb{Z}$ appears either $\lfloor v/M \rfloor$ or $\lceil v/M \rceil$ times.

\begin{definition}
Let $p$ be an odd prime, let $d$ be a positive integer, and let $q = p^d.$ Let $\alpha$ be a primitive element of $\mathbb{F}_q,$ and let $M|q-1.$ Following \cite{G5}, we set 
\begin{eqnarray*}
D_k = \lbrace \alpha^{Mi + k}-1\mid 0 \leq i < (q-1)/M \rbrace \mbox{ ~for~ }  0 \leq k \leq M-1.
\end{eqnarray*}
An M-ary Sidelnikov sequence $\mathbf{s} = s_0 s_1 s_2 \ldots$ is a sequence of period $q-1$ whose first $q-1$ elements are defined as follows: for $0 \leq j < q-1$
\[s_j = \begin{cases} 0 & \text{ if } \alpha^j = -1,\\ k & \text{ if } \alpha^j \in D_k. \end{cases}\] 
\end{definition}

This class of sequences was originally discovered by Sidelnikov in 1969 \cite{S2} and, in the binary case (i.e. the case in which $M = 2$) rediscovered independently by Lempel, Cohn, and Eastman in 1977 \cite{C3}.

It is known that the Sidelnikov sequences have low out-of-phase autocorrelation. Indeed, their out-of-phase autocorrelation values have magnitude at most 4 (see, for instance, \cite{C1}). It is clear from the definition of these sequences that they have the balance property.

\begin{remark}
If $c \in \mathbb{Z}/M\mathbb{Z},$ then we stipulate that $c\mathbf{a}$ is the sequence whose $i$th entry is $ca_i$ and we say that $c\mathbf{a}$ is a \emph{constant multiple} of $\mathbf{a}$.
\end{remark} 

The authors of \cite{K4} consider the family $\mathcal{S}$ of sequences consisting of all nonzero constant multiples of an $M$-ary Sidelnikov sequence. They are able to prove that $\mathcal{S}$ has rather remarkable correlation properties.  Indeed, they use the Weil bound to prove the following upper bound on the cross-correlation of a two distinct constant multiples of a Sidelnikov sequence.
\begin{theorem} \label{good bound} \cite[Theorem 3]{K4} Let $q$ be a power of an odd prime, and let $M|q-1.$ Let $\mathbf{s}$ be an M-ary Sidelnikov sequence over $\mathbb{F}_q$. Let $c_1, c_2 \in \mathbb{Z}/M\mathbb{Z}$, $c_1\not=c_2$ and $c_1, c_2 \neq 0$. Let $\mathbf{a}:= c_1\mathbf{s}$, and let $\mathbf{b}:= c_2 \mathbf{s}$. Then for each $\tau = 0, \ldots ,q-2,$ 
\[|\mathcal{C}_{\mathbf{a},\mathbf{b}}(\tau)| \leq \sqrt{q} + 3.\]
\end{theorem}

As the authors of \cite{K4} note, the upper bound in Theorem 1.2 implies that $\mathcal{S}$ is a family of $M-1$ sequences with (asymptotically) nearly optimal cross-correlation. For, if $\mathcal{C}_{\text{max}}$ denotes the maximum of the cross-correlation and out-of-phase autocorrelation values of the sequences in $\mathcal{S},$ then by Theorem \ref{welch}, 
\[\mathcal{C}_{\text{max}} \geq \sqrt{\frac{(q-1)^2(M-2)}{(q-1)(M-1)-1}},\]
and this expression is asymptotically equivalent to $\sqrt{q}.$\\

But the suitability of a sequence family for CDMA depends not just on the maximum cross-correlation value of its sequence pairs (although the maximum value is important). It is also important to determine the rest of the cross-correlation values of the pairs from the sequence family. Furthermore, we note that for fixed values of $M$ and $q,$ there is some distance between the lower bound given in Theorem \ref{welch} and the upper bound given in Theorem \ref{good bound}. For instance, if $q = 5^5$ and $M = 4,$ then 
\[\sqrt{\frac{(q-1)^2(M-2)}{(q-1)(M-1)-1}} \approx 45.63 \mbox{~ and~ } \sqrt{q}+3 \approx 58.90.\] 
So, the known bounds given do not completely specify the maximum cross-correlation values. 
Thus, it is of interest to determine the precise cross-correlation values of the sequences from $\mathcal{S}$. 

In this paper, we show that the cross-correlation values of two sequences from $\mathcal{S}$ can be represented in terms of certain Jacobi sums. 
We then use known evaluations of Jacobi sums to determine the precise cross-correlation distribution of $\mathcal{S}$ in certain cases. 
  
\section{Preliminaries}
Let $q$ be a prime power, and let $\mathbb{F}_q$ be a finite field with $q$ elements. A homomorphism $\chi:\mathbb{F}_q^* \to \mathbb{C}$ is called a \emph{character} on $\mathbb{F}_q^*$. The trivial character on $\mathbb{F}_q^*$ is defined by the rule that $\chi(x)=1$ for all $x\in\mathbb{F}_q^*$. Any character on $\mathbb{F}_q^*$ which is different from the trivial character is called a \emph{nontrivial character}. 

We note that if $\chi$ is a nontrivial character on $\mathbb{F}_q^*$, it is common to extend $\chi$ to a map on $\mathbb{F}_q$ by setting $\chi(0) = 0$. However, it is sometimes useful to define $\chi(0)$ to be equal to something else. 
In this paper, we consider both characters $\chi$ on $\mathbb{F}_q$ for which $\chi(0) = 0$ and characters $\chi$ on $\mathbb{F}_q$ for which $\chi(0) = 1$. 

It is possible to define a logarithm over $\mathbb{F}_q.$ Let $\alpha$ be a primitive element of $\mathbb{F}_q.$ For $x \in \mathbb{F}_q,$ we stipulate that 
\[\text{log}_{\alpha}(x) = \begin{cases} i & \text{if } x = \alpha^i, \hspace{0.05in} 0 \leq i \leq q-2, \\
0 & \text{if } x = 0.\end{cases}\]
As Gong and Yu note (\cite{G5} and \cite{G7}) the $M$-ary Sidelnikov sequence $\mathbf{s}$ defined over $\mathbb{F}_q$ using $\alpha$ is completely determined by the congruences 
\begin{equation} 
\label{id1} s_i \equiv \text{log}_{\alpha}(\alpha^i+1) \pmod{M}, \hspace{1cm} 0 \leq i \leq q-2.
\end{equation}
Still following Gong and Yu, we define a multiplicative character $\psi_M$ of order $M$ on $\mathbb{F}_q$ by the rule that for $x \in \mathbb{F}_q$, 
\begin{equation} 
\label{id2} \psi_M(x) = \text{exp}\Big( \frac{2\pi i\text{log}_{\alpha}(x)}{M}\Big).
\end{equation}
Note that $\psi_M(0) = 1$.
Gong and Yu remark that (\ref{id1}) and (\ref{id2}) imply the identity
\begin{equation} 
\label{id3} \text{exp}\Big(\frac{2\pi i s_j}{M}\Big) = \psi_M(\alpha^j + 1), \hspace{0.1in} 0 \leq j \leq q-2. 
\end{equation}
Let $\chi_M$ denote the character on $\mathbb{F}_q$ defined by
\begin{equation*} 
\chi_M(x) = \begin{cases} \psi_M(x) & \text{if } x \not= 0, \\ 0 &  \text{if } x = 0.\end{cases}
\end{equation*} 

We shall need to make use of certain types of character sums called \emph{Gauss sums} and \emph{Jacobi sums}. These sums were originally studied by C. F. Gauss and C. G. I. Jacobi, hence the names. 

\begin{definition}
Let $q = p^k$ be a power of a prime $p$, and let $\chi$ and $\psi$ be nontrivial characters on $\mathbb{F}_q$ {\rm (}which map $0$ to $0${\rm )}. Then the sum 
\[J(\chi,\psi) = \sum_{x \in \mathbb{F}_q} \chi(x)\psi(1-x)\]
is called a \emph{Jacobi sum}, and the sum 
\[G(\chi) = \sum_{x \in \mathbb{F}_q} \chi(x)\text{exp}\Big(\frac{2\pi i \text{Tr}(x)}{p}\Big)\]
 {\rm (}where Tr denotes the field trace from $\mathbb{F}_q$ to $\mathbb{F}_p${\rm)} is called a \emph{Gauss sum}.
\end{definition}

The following identity relates Gauss and Jacobi sums, see \cite[Theorem 2.1.3]{I1}. If $\chi\psi$ is not the trivial character, then 
\begin{equation} 
\label{Gauss Jacobi} J(\chi,\psi) = \frac{G(\chi)G(\psi)}{G(\chi\psi)}. 
\end{equation}
We note the following basic properties of Gauss and Jacobi sums (see \cite{I1}).
\begin{eqnarray*}
&&G(\chi) = q-1 \mbox{ if } \chi \mbox{ is trivial}.\\ 
&&|G(\chi)| = \sqrt{q} \mbox{ if } \chi  \mbox{ is nontrivial}.\\
&&J(\chi,\psi) = -1 \mbox{ if } \chi = \overline{\psi}.\\
&&|J(\chi,\psi)| = \sqrt{q} \mbox{ if } \chi \neq \overline{\psi}.
\end{eqnarray*}

It is not known, in general, how to produce exact, simple, closed form evaluations of $J(\chi,\psi)$ and $G(\chi).$ That being said, the sums $J(\chi,\psi)$ and $G(\chi)$ have been evaluated in certain special cases.
One case in which there are known evaluations of Gauss and Jacobi sums is the small order case (when the order of the characters involved in the sums is small relative to $q$). Explicit formulae for Gauss and Jacobi sums over characters of orders $3,5,7,\ldots$ etc have been given in the literature (see \cite[Chapters 3 and 4]{B1}). However, these formulae are rather lengthy and difficult to apply. Alternatively, there are algorithms that can be used to compute Jacobi sums over characters of small order (see, for instance, \cite{B2} and \cite{V1}). Let $\chi$ and $\psi$ be two characters on $\mathbb{F}_q,$ and let $e = \text{gcd}(|\chi|,|\psi|).$ The algorithm given in \cite{V1} computes the Jacobi sum $J(\chi,\psi)$ faster than just naively summing the series so long as $\phi(e)^{\phi(e)}\leq q$, where $\phi$ is the Euler totient function. There are other situations in which it is possible to give explicit evaluations of Gauss and Jacobi sums. We will make use of known evaluations for the so-called pure Gauss and Jacobi sums. A Gauss or Jacobi sum is called \emph{pure} if some positive integral power of it is real. The following theorem completely classifies pure Gauss sums (see \cite[Theorems  11.6.3 and 11.6.4]{B1}). 

\begin{theorem} \label{pureGauss} Let $m \in \mathbb{N},$ let $p$ be an odd prime, and let $q = p^m.$ Let $k\mid q-1$, and let $\chi$ be a character of order $k$. Then $G(\chi)$ is pure if and only if there exists a positive integer $x$ such that $p^x \equiv -1 \pmod{k}.$ Furthermore, if there exist such integers and $t$ is the least such integer, then there exists a positive integer $s$ such that $m = 2ts,$ and
\[G(\chi) =  (-1)^{s-1+(p^t+1)s/k}p^{m/2}.\]
\end{theorem}
We can use Theorem \ref{pureGauss} to evaluate a large class of pure Jacobi sums. 

\begin{corollary} \label{longCor} Let $m \in \mathbb{N},$ let $p$ be an odd prime, and let $q = p^m$. Let $k\mid q-1$, and let $\chi$ be a character of order $k$. Suppose that there exists a positive integer $x$ such that $p^x \equiv -1 \pmod{k}$. Let $a$ and $b$ be integers such that $a+b \not\equiv 0 \pmod{k}$. Then there exist positive integers $w$, $y$, and $z$ such that $p^w \equiv -1 \pmod{k/\text{gcd}(k,a)}$, $p^y \equiv -1 \pmod{k/\text{gcd}(k,b)}$, and $p^z \equiv -1 \pmod{k/\text{gcd}(k,a+b)}$. Let $t_a$, $t_b$, and $t_{a+b}$ be the least such integers. Then there exist positive integers $s_a$, $s_b$, and $s_{a+b}$ such that $m = 2t_as_a = 2t_bs_b = 2t_{a+b}s_{a+b}$. Furthermore, 
\[J(\chi^a,\chi^b) = (-1)^{s_a + s_b + s_{a+b} + 1 + ((p^{t_a}+1)s_a\text{gcd}(k,a) + (p^{t_b}+1)s_b\text{gcd}(k,b) + (p^{t_{a+b}}+1)s_{a+b}\text{gcd}(k,a+b))/k}p^{m/2}.\]
\end{corollary}

\emph{Proof.}
For each positive integer $c,$ $|\chi^c| = k/\text{gcd}(k,c)$ divides $k.$ Hence, there exists a positive integer $x$ such that $p^x \equiv -1 \pmod{k/\text{gcd}(k,c)}.$ Let $t_c$ be the least such integer. Then by Theorem \ref{pureGauss}, there exists an integer $s_c$ such that $m = 2t_cs_c$ and 
\begin{eqnarray}
G(\chi^c) = (-1)^{s_c-1+(p^{t_c}+1)s_c\text{gcd}(k,c)/k}p^{m/2}.
\end{eqnarray}
By (2.4) and (2.5) we have
\begin{eqnarray*}
&&\hspace{-1cm}J(\chi^a,\chi^b)
= \frac{G(\chi^a)G(\chi^b)}{G(\chi^{a+b})}\\
&=& \frac{(-1)^{s_a-1+(p^{t_a}+1)s_a\text{gcd}(k,a)/k}p^{m/2}(-1)^{s_b-1+(p^{t_b}+1)s_b\text{gcd}(k,b)/k}p^{m/2}}{(-1)^{s_{a+b}-1+(p^{t_{a+b}}+1)s_{a+b}\text{gcd}(k,a+b)/k}p^{m/2}}\\
&=& (-1)^{s_a + s_b + s_{a+b} + 1 + ((p^{t_a}+1)s_a\text{gcd}(k,a) + (p^{t_b}+1)s_b\text{gcd}(k,b) + (p^{t_{a+b}}+1)s_{a+b}\text{gcd}(k,a+b))/k}p^{m/2}. \qed
\end{eqnarray*}

The Jacobi sum evaluation given in Corollary \ref{longCor} is a bit cumbersome, but in some cases it simplifies nicely (as in the following corollary).

\begin{corollary} \label{simpleCor} Let $m \in \mathbb{N}$ be such that $m \equiv 0 \pmod{4}$, let $p$ be an odd prime, and let $q = p^m$. Let $k\mid q-1$, and let $\chi$ be a character of order $k$. Suppose that there exists a positive integer $x$ such that $p^x \equiv -1 \pmod{k}$. Furthermore, suppose that the least positive integer satisfying this equation is odd. Let $a$ and $b$ be integers such that $a+b \not\equiv 0 \pmod{k}$. Then 
\[J(\chi^a,\chi^b) = -p^{m/2}.\]
\end{corollary}

{\em Proof.} Let $c$ be a positive integer, and suppose that $t$ is the least positive integer satisfying $p^x \equiv -1 \pmod{k}.$ As above, let $t_c$ be the least positive solution of $p^y \equiv -1 \pmod{k/\text{gcd}(k,c)}.$ Note that, in this case, the order of $p$ in the multiplicative group $(\mathbb{Z}/\text{gcd}(k,c)\mathbb{Z})^*$ is $2t_c.$ By the Division Algorithm, there exist unique integers $y$ and $r$ such that $t = t_cy + r$ and $0 \leq r <t_c$. Thus, 
\begin{eqnarray}
-1 \equiv p^t \equiv p^{t_cy+r} \equiv (-1)^yp^r \pmod{k/\text{gcd}(k,c)}.
\end{eqnarray}
If $y$ is even, then $p^r \equiv -1 \pmod{k/\text{gcd}(k,c)}.$ But since $r<t_c,$ this contradicts the fact that $t_c$ is the least positive integer satisfying $p^y \equiv -1 \pmod{k/\text{gcd}(k,c)}$. Hence $y$ must be odd. Then by (2.6) we have  $p^r \equiv 1 \pmod{k/\text{gcd}(k,c)}$. Thus, unless $r = 0,$ since $r<t_c<2t_c,$ this contradicts the fact that $2t_c$ is the order of $p$ in $(\mathbb{Z}/\text{gcd}(k,c)\mathbb{Z})^*.$ It follows that $r = 0.$ Hence, $t_c|t.$
Thus, if $t$ is odd, then the terms $t_a,$ $t_b,$ and $t_{a+b}$ in Corollary \ref{longCor} must also be odd. Hence, if $m \equiv 0 \pmod{4},$ then $s_a,$ $s_b,$ and $s_{a+b}$ must all be even. The assertion now follows from Corollary 2.1.
\qed \\

We also need to make use of cyclotomic numbers. 

\begin{definition}
Let $q = p^m$ be an odd prime power, let $k\mid q-1$, and let $\alpha$ be a primitive element in $\mathbb{F}_q$. Then for $0 \leq u \leq k-1$, we set 
\[C_u = \left\{ \alpha^{k\ell + u} \mid 0 \leq \ell \leq \frac{q-1}{k} \right\}.\] 
The sets $C_u$ are called \emph{cyclotomic classes}. For $0 \leq u,v \leq k-1$, the \emph{cyclotomic number} $(u,v)_k$ is defined to be the number of elements $x \in C_u$ such that $1+x \in C_v$. 
\end{definition}

Cyclotomic numbers were first studied by C. F. Gauss as part of his work on the problem of finding straight edge and compass constructions for regular $n$-gons \cite{S3}.
There is a connection between cyclotomic numbers and Jacobi sums. For a proof of the next result, see \cite[Theorem 2.5.1]{B1}.
\begin{theorem} \label{cyclotomic Jacobi}
Let $q = p^m$ be an odd prime power, let $k\mid q-1$, and let $\alpha$ be a primitive element in $\mathbb{F}_q$. Let $\chi$ be a character of order $k$. Then for $0 \leq u \leq k-1$, 
\[k^2(u,v)_k = \sum_{w,x \pmod{k}} \chi^u(-1)(J(\chi^u,\chi^v)-\ell)\zeta^{-uw -vx},\]
where $\zeta$ is a certain primitive $k$th root of unity and $\ell$ is the number of trivial characters in the set $\lbrace \chi^u,\chi^v \rbrace.$
\end{theorem}
A number of authors have made use of Theorem \ref{cyclotomic Jacobi} to compute formulae for cyclotomic numbers of ``small order.'' See \cite{S3} for a summary of many such results. By making use of the evaluations of Jacobi sums in the pure case, one can deduce the following corollary to Theorem \ref{cyclotomic Jacobi} (see \cite[Section 11.6]{B1} and \cite{E1}).
\begin{corollary} \label{cyc numbers} Let $m \in \mathbb{N}$, let $p$ be an odd prime, and let $q = p^m$. Let $k\mid q-1$, and suppose that there exists a positive integer $x$ such that $p^x \equiv -1 \pmod{k}$. Let $t$ be the least such integer. Then there exists a positive integer $s$ such that $m = 2ts$. Furthermore, 
\[k^2(0,0)_k = q+1-3k-(k-1)(k-2)(-1)^sp^{m/2}.\]
If $b \not\equiv 0 \pmod{k}$, then
\[k^2(0,b)_k = k^2(b,0)_k = k^2(k-b,k-b)_k = q+1-k+(-1)^s(k-2)p^{m/2}.\]
If $a,b,a-b \not\equiv 0 \pmod{k}$, then 
\[k^2(a,b)_k = q+1-2(-1)^sp^{m/2}.\]
\end{corollary}

\section{Cross-correlation distributions}
\begin{theorem} \label{corr} Let $q$ be an odd prime power, let $\alpha$ be a primitive element of $\mathbb{F}_q$, and let $M\mid q-1$. Let $c_1,c_2 \in \lbrace 1,\ldots,M-1 \rbrace$ with $c_1 \neq c_2$. Let $\mathbf{s}$ be the $M$-ary Sidelnikov sequence over $\mathbb{F}_q$ defined using $\alpha.$ Then 
\[\mathcal{C}_{c_1\mathbf{s},c_2\mathbf{s}}(0) = -1\]
and for $\tau = 1,\ldots,q-2$, 
\[\mathcal{C}_{c_1\mathbf{s},c_2\mathbf{s}}(\tau) = \chi_M^{c_1}(1-\alpha^{-\tau})\chi_M^{-c_2}(1-\alpha^{\tau})(J(\chi_M^{c_1},\chi_M^{-c_2})+2)-1.\]
\end{theorem}

\emph{Proof.} For any $\tau$, by (1.1) and (\ref{id3}) we have
\begin{eqnarray*}
\mathcal{C}_{c_1\mathbf{s},c_2\mathbf{s}}(\tau) 
= \sum_{t = 0}^{q-2} \text{exp}\Big({\frac{2\pi i(c_1s_t - c_2s_{t + \tau})}{M}}\Big) 
= \sum_{t=0}^{q-2} \psi_M^{c_1}(\alpha^t+1) \psi_M^{-c_2}(\alpha^{t+\tau}+1)\\
\end{eqnarray*}
As $t$ ranges from $0$ to $q-2$, $\alpha^t +1$ ranges over all of the elements of $\mathbb{F}_q$ except for $1$. So, 
\begin{eqnarray*}
\mathcal{C}_{c_1\mathbf{s},c_2\mathbf{s}}(0) = \sum_{t = 0}^{q-2} \psi_M^{c_1-c_2}(\alpha^t + 1) = -1.
\end{eqnarray*}
Now, assume that $\tau = 1,...,q-2.$ Note that 
\begin{eqnarray*}
\mathcal{C}_{c_1\mathbf{s},c_2\mathbf{s}}(\tau) &=& \psi_M^{-c_2}(\alpha^{\tau})\sum_{t=0}^{q-2} \psi_M^{c_1} (\alpha^t+1)\psi_M^{-c_2}(\alpha^t + 1 + (\alpha^{-\tau} -1)).
\end{eqnarray*}
Again, using the fact that as $t$ ranges from $0$ to $q-2$, $\alpha^t +1$ ranges over all of the elements of $\mathbb{F}_q$ except for $1,$ we get that 
\begin{eqnarray*}
\mathcal{C}_{c_1\mathbf{s},c_2\mathbf{s}}(\tau) 
&=& \psi_M^{-c_2}(\alpha^{\tau})\big( \sum_{x \in \mathbb{F}_q} \psi_M^{c_1} (x)\psi_M^{-c_2}(x + (\alpha^{-\tau} -1)) - \psi_M^{-c_2}(\alpha^{-\tau})\big)\\
&=& \psi_M^{-c_2}(\alpha^{\tau})\sum_{x \in \mathbb{F}_q} \psi_M^{c_1} (x)\psi_M^{-c_2}(x + (\alpha^{-\tau} -1)) - 1.
\end{eqnarray*}
Finally, as $x$ runs over the elements of $\mathbb{F}_q$, $-(\alpha^{-\tau}-1)x$ also runs over the elements of $\mathbb{F}_q$. So,
\begin{eqnarray*} 
\mathcal{C}_{c_1\mathbf{s}, c_2\mathbf{s}}(\tau) 
&=&  \psi_M^{-c_2}(\alpha^{\tau})\sum_{x \in \mathbb{F}_q} \psi_M^{c_1} (-(\alpha^{-\tau}-1)x)\psi_M^{-c_2}(-(\alpha^{-\tau}-1)x + (\alpha^{-\tau} -1)) - 1 \\
&=& \psi_M^{-c_2}(\alpha^{\tau})\psi_M^{c_1}(-(\alpha^{-\tau}-1)) \psi_M^{-c_2}(\alpha^{-\tau}-1)\sum_{x \in \mathbb{F}_q} \psi_M^{c_1}(x)\psi_M^{-c_2}(1-x)-1 \\
&=& \chi_M^{c_1}(1-\alpha^{-\tau})\chi_M^{-c_2}(1-\alpha^{\tau})(J(\chi_M^{c_1},\chi_M^{-c_2})+2)-1. \qed
\end{eqnarray*}
\vspace{2mm}

So long as $\phi(M)^{\phi(M)} \leq q,$ Theorem \ref{corr} can be used in concert with the algorithm from \cite{V1} to facilitate computations of cross-correlations of constant multiples of M-ary Sidelnikov sequences. Indeed, we consider this to be the main application of Theorem \ref{corr}. 

We will also use Theorem \ref{corr} together with known evaluations of pure Jacobi sums to deduce explicit formulae for the cross-correlation values of constant multiples of M-ary Sidelnikov sequences in certain special cases.
Note that  
\begin{equation} 
\label{conj} J(\chi_M^{c_1},\chi_M^{-c_2}) = -1 \mbox{ if } c_1 + c_2 \equiv 0 \pmod{M}.
\end{equation}
The next corollary follows directly from (\ref{conj}), Theorem \ref{corr}, Corollary \ref{longCor}, and Corollary \ref{simpleCor}. 
\begin{corollary} \label{longerCor} Let $m \in \mathbb{N}$, $p$ be an odd prime, and $q = p^m$. Let $\alpha$ be a primitive element of $\mathbb{F}_q,$ and let $M\mid q-1$. Let $c_1,c_2 \in \lbrace 1,\ldots,M-1 \rbrace$ with $c_1 \neq c_2$. Let $\mathbf{s}$ be the $M$-ary Sidelnikov sequence over $\mathbb{F}_q$ defined using $\alpha$. Let $\tau \in\{1,\ldots,q-2\}$.\\
\noindent If $c_1+ c_2 \equiv 0 \pmod{M}$, then 
\[\mathcal{C}_{c_1\mathbf{s},c_2\mathbf{s}}(\tau) = \chi_M^{c_1}(1-\alpha^{-\tau})\chi_M^{-c_2}(1-\alpha^{\tau})-1.\]
\noindent Suppose that $c_1 + c_2 \not\equiv 0 \pmod{M}.$ Suppose also that there exists a positive integer $x$ such that $p^x \equiv -1 \pmod{M}.$ Then there exist positive integers $w,$ $y,$ and $z$ such that $p^w \equiv -1 \pmod{M/\text{gcd}(M,c_1)},$ $p^y \equiv -1 \pmod{M/\text{gcd}(M,c_2)},$ and $p^z \equiv -1 \pmod{M/\text{gcd}(M,c_1+c_2)}.$ Let $t_{c_1},$ $t_{c_2},$ and $t_{c_1+c_2}$ be the least such integers. Then there exist positive integers $s_{c_1},$ $s_{c_2},$ and $s_{c_1+c_2}$ such that 
\begin{eqnarray*}
m = 2t_{c_1}s_{c_1} = 2t_{c_2}s_{c_2} = 2t_{c_1+c_2}s_{c_1+c_2}.
\end{eqnarray*}
 Furthermore, 
\begin{eqnarray*}
\mathcal{C}_{c_1\mathbf{s},c_2\mathbf{s}}(\tau) = \chi_M^{c_1}(1-\alpha^{-\tau})\chi_M^{-c_2}(1-\alpha^{\tau}) ((-1)^{\epsilon}p^{m/2} + 2)-1,
\end{eqnarray*}
where
\begin{eqnarray*}
\epsilon&=&s_{c_1} + s_{c_2} + s_{c_1+c_2}  + 1 + ((p^{t_{c_1}}+1)s_{c_1}\text{gcd}(M,c_1) + (p^{t_{c_2}}+1)s_{c_2}\text{gcd}(M,c_2)\\
&&+ (p^{t_{c_1+c_2}}+1)s_{c_1+c_2}\text{gcd}(M,c_1+c_2))/M.
\end{eqnarray*}
If $m \equiv 0 \pmod{4}$ and the least positive integer satisfying the equation $p^x \equiv -1 \pmod{M}$ is odd, then 
\[\mathcal{C}_{c_1\mathbf{s},c_2\mathbf{s}}(\tau) = \chi_M^{c_1}(1-\alpha^{-\tau})\chi_M^{-c_2}(1-\alpha^{\tau})(-p^{m/2} + 2) -1.\]
\end{corollary}
By borrowing some ideas from the authors of \cite{C1}, we can use Corollary \ref{longerCor} to determine the cross-correlation distributions of certain families of constant multiples of Sidelnikov sequences. 

By Corollary 3.1, for fixed $c_1$ and $c_2$,  the cross-correlation $C_{c_1\mathbf{s},c_2\mathbf{s}}(\tau)$ depends on $\chi_M^{c_1}(1-\alpha^{-\tau})\chi_M^{-c_2}(1-\alpha^{\tau}).$ 
For $\tau \in\{1,2,\ldots,q-2\}$, let $y = \alpha^{\tau}$. Note that $y \neq 0,1$. We set 
\begin{equation}
\omega_M = \text{exp}\big(\frac{2\pi i}{M}\big).
\end{equation}
Let $0 \leq u,v \leq M-1$ be such that  
\begin{equation}
\chi_M^{c_1}\Big(\frac{y-1}{y}\Big) = \omega_M^u \hspace{0.1in} \text{and} \hspace{0.1in} \chi_M^{c_2}\Big(\frac{1}{1-y}\Big) = \omega_M^v.
\end{equation}
From (3.3) we have
\begin{eqnarray}
\chi_M^{c_1}(1-\alpha^{-\tau})\chi_M^{-c_2}(1-\alpha^{\tau}) &=& \chi_M^{c_1}(\alpha^{-\tau}(\alpha^{\tau}-1))\chi_M^{-c_2}(1-\alpha^{\tau})\\
&=& \chi_M^{c_1}\Big(\frac{y-1}{y}\Big)\chi_M^{c_2}\Big(\frac{1}{1-y}\Big) \nonumber \\
&=& \omega_M^{c_1u + c_2v}.\nonumber
\end{eqnarray}
Hence, the values taken on by $\chi_M^{c_1}(1-\alpha^{-\tau})\chi_M^{-c_2}(1-\alpha^{\tau})$ are completely determined by the cardinalities of the sets 
\[S_{u,v} = \left\{ y \in \mathbb{F}_{p^m}\setminus \lbrace 0,1 \rbrace \mid \chi_M\Big(\frac{y-1}{y}\Big) = \omega_M^u \hspace{0.1in} \text{and} \hspace{0.1in} \chi_M\Big(\frac{1}{1-y}\Big) = \omega_M^v \right\},\]
where $0 \leq u,v \leq M-1.$ The authors of \cite{C1} characterized the cardinalities of these sets in terms of cyclotomic numbers.
\begin{theorem} \label{cyc char} \cite[Theorem 11]{C1} Let $0 \leq u,v \leq M-1$. Then 
\[|S_{u,v}| = (u+v,v)_M.\]
\end{theorem}
In the following example we use  Corollary \ref{longerCor}, Theorem \ref{cyc char}, and Corollary \ref{cyc numbers} to illustrate how to determine the cross-correlation distribution of the family of constant multiples of an $M$-ary Sidelnikov sequence.

\begin{example} Let $p = 3$, $m = 4$, $q = p^m = 81$, and $M = 4$. 
By {\rm(3.2)}, $\omega_M = i$. Let $\alpha$ be a primitive element of $\mathbb{F}_{81}$, and let $\mathbf{s}$ be the $M$-ary Sidelnikov sequence defined using $\alpha$. Note that $M\mid q-1$. Furthermore, $p \equiv -1 \pmod{M}$, so that $t = 1$. We also have that  $m = 2t \cdot 2$, so that $s = 2$. Thus, we can apply {\rm Corollary \ref{longerCor}} to determine the cross-correlation values of any two constant multiples of $\mathbf{s}$. 
We will explicitly show how this is done for $c_1 = 1$ and $c_2 = 2$.

Since $c_1 + c_2 \not\equiv 0 \pmod{M},$ for $\tau\in\lbrace 1,2,\ldots,q-2 \rbrace$, by {\rm Corollary 3.1} we have 
\begin{equation} 
C_{\mathbf{s},2\mathbf{s}}(\tau) =  \chi_M^{c_1}(1-\alpha^{-\tau})\chi_M^{-c_2}(1-\alpha^{\tau})(-9+2) -1.
\end{equation}
For $0 \leq u,v \leq 3$, by {\rm(3.4)} and {\rm Theorem 3.2} we have
\begin{equation} 
\chi_M^{c_1}(1-\alpha^{-\tau})\chi_M^{-c_2}(1-\alpha^{\tau}) = \omega_M^{u+2v} = i^{u+2v}
\end{equation}
for exactly $(u+v,v)_M$ values of $\tau$. Note that
\begin{equation*} 
u + 2v \equiv 0 \pmod{4} \mbox{ if and only if } 
\begin{cases} u = v = 0,  \mbox{ or} \\ u = 0 \mbox{ and } v = 2,  \mbox{ or} \\  u = 2 \mbox{ and } v = 1,  \mbox{ or}\\ u = 2 \mbox{ and } v = 3.
\end{cases}
\end{equation*} 
By {\rm Corollary \ref{cyc numbers}}, we obtain
\begin{equation}
(0,0)_4 + (2,2)_4 + (3,1)_4 + (1,3)_4 = 15.
\end{equation}
Thus by {\rm(3.6)} and {\rm(3.7)}, we have
\[\chi_M^{c_1}(1-\alpha^{-\tau})\chi_M^{-c_2}(1-\alpha^{\tau}) = 1\]
for exactly $15$ values of $\tau$. \\

Note that
\begin{equation*} 
u + 2v \equiv 1 \pmod{4} \mbox{ if and only if } 
\begin{cases} u =1   \mbox{ and } v = 0, \mbox{ or} \\ u = 1 \mbox{ and } v = 2, \mbox{ or} \\  u = 3 \mbox{ and } v = 1, \mbox{ or} \\ u = 3 \mbox{ and } v = 3.
\end{cases}
\end{equation*} 
By {\rm Corollary \ref{cyc numbers}}, we obtain
\begin{eqnarray}
(1,0)_4 + (3,2)_4 + (0,1)_4 + (2,3)_4 = 20.
\end{eqnarray}
Thus by {\rm(3.6)} and {\rm(3.8)}, we have
\[\chi_M^{c_1}(1-\alpha^{-\tau})\chi_M^{-c_2}(1-\alpha^{\tau}) = i\]
for exactly $20$ values of $\tau$. \\

Note that
\begin{equation*} 
u + 2v \equiv 2 \pmod{4} \mbox{ if and only if } 
\begin{cases} u =0   \mbox{ and } v = 1, \mbox{ or} \\ u = 0 \mbox{ and } v = 3, \mbox{ or} \\  u = 2 \mbox{ and } v = 0, \mbox{ or} \\ u = 2 \mbox{ and } v = 2.
\end{cases}
\end{equation*} 
By {\rm Corollary \ref{cyc numbers}}, we obtain
\begin{eqnarray}
(1,1)_4 + (3,3)_4 + (2,0)_4 + (0,2)_4 = 24.
\end{eqnarray}
Thus by {\rm(3.6)} and {\rm(3.9)}, we have
\[\chi_M^{c_1}(1-\alpha^{-\tau})\chi_M^{-c_2}(1-\alpha^{\tau}) = -1\]
for exactly $24$ values of $\tau$.\\ 

Note that
\begin{equation*} 
u + 2v \equiv 3 \pmod{4} \mbox{ if and only if } 
\begin{cases} u = v = 1, \mbox{ or} \\ u = 1 \mbox{ and } v = 3, \mbox{ or} \\  u = 3 \mbox{ and } v = 0, \mbox{ or} \\ u = 3 \mbox{ and } v = 2.
\end{cases}
\end{equation*} 
By {\rm Corollary \ref{cyc numbers}}, we obtain
\begin{eqnarray}
(2,1)_4 + (0,3)_4 + (3,0)_4 + (1,2)_4 = 20.
\end{eqnarray}
Thus by {\rm(3.6)} and {\rm(3.10)}, we have 
\[\chi_M^{c_1}(1-\alpha^{-\tau})\chi_M^{-c_2}(1-\alpha^{\tau}) = -i\]
for exactly $20$ values of $\tau$. 

Thus, as $\tau$ ranges over $\lbrace 1, 2,\ldots,79\rbrace$, by {\rm(3.5)} we have
\[C_{\mathbf{s},2\mathbf{s}}(\tau) =
 \begin{cases} 
 -8 & \text{exactly $15$ times} \\
-7i-1 & \text{exactly $20$ times} \\
6 & \text{exactly $24$ times} \\
7i-1 & \text{exactly $20$ times}.
\end{cases} 
\]
We also note that
\[C_{\mathbf{2s},3\mathbf{s}}(\tau)=C_{\mathbf{s},2\mathbf{s}}(\tau)\]
and
\[C_{\mathbf{s},3\mathbf{s}}(\tau) =
 \begin{cases} 
 0 & \text{exactly $15$ times} \\
i-1 & \text{exactly $20$ times} \\
-2 & \text{exactly $24$ times} \\
-i-1 & \text{exactly $20$ times}.
\end{cases} 
\]

\end{example}

\section*{Acknowledgements} 
The research of Ay\c{s}e Alaca was supported by a Discovery Grant
from the Natural Sciences and Engineering Research Council of Canada (RGPIN-418029-2013).

\end{document}